\hoffset=1.5cm
\newcount\refco
\refco=0

\advance\refco by 1
\newcount\juana\juana=\refco

\advance\refco by 1
\newcount\baxter\baxter=\refco

\advance\refco by 1
\newcount\koles\koles=\refco

\advance\refco by 1
\newcount\lodayuno\lodayuno=\refco

\advance\refco by 1
\newcount\lodaydos\lodaydos=\refco

\advance\refco by 1
\newcount\operads\operads=\refco

\advance\refco by 1
\newcount\sib\sib=\refco

\newcount\thco
\thco 0\relax
\newcount\deco
\deco 0\relax
\newcount\remco
\remco 0\relax
\newcount\cf
\cf 0\relax
\newcount\reone
\reone=0
\newcount\kaka
\kaka=0
\newcount\sectno
\sectno=0
\newcount\coco
\coco=0
\newcount\prco
\prco=0

\def\ponum{\advance \cf by 1
\eqno{(\number\cf)}}

\def\iz{\dashv}
\def\de{\vdash}
\def\th#1{
\advance\thco 1
{\bf Theorem \number\thco.} {\it #1}
}

\def\pr#1{
\advance\prco 1
{\bf Proposition \number\prco.} {\it #1}
}

\def\defi#1{
\advance\deco 1
{\bf Definition \number\deco.} {\it #1}}

\def\rem#1{
\advance\remco 1
{\bf Remark \number\remco.} { #1}
}

\def\co#1{
\advance\coco 1
{\bf Corollary \number\coco.} {\it #1}
}

\def\raniz{\mathop{\hbox{Rann}^{\iz}}}
\def\rande{\mathop{\hbox{Rann}^{\de}}}
\def\laniz{\mathop{\hbox{Lann}^{\iz}}}
\def\lande{\mathop{\hbox{Lann}^{\de}}}
\def\ann{\mathop{\hbox{Ann}}}

\def\op#1{{#1}^{\hbox{\sevenrm op}}}

\hsize=13cm
\centerline{\bf Associative dialgebras from a structural viewpoint}

\footnote{}{\sevenrm The authors have been supported by the Spanish MEC and Fondos
FEDER through projects MTM\-2007\-60333 and
MTM2010\-15223, jointly by the Junta de Andaluc\'{\i}a and Fondos FEDER
through projects FQM-336, FQM-2467 and
FQM-3737 and by the Spanish Ministry of Education and Science under
project ``Ingenio Mathematica (i-math)''
No. CSD2006-00032 (Consolider-Ingenio 2010).}

\bigskip
\centerline{\sl C\'andido Mart\'{\i}n Gonz\'alez}
\bigskip
\centerline{\sevenrm Department of Algebra, Geometry and Topology}
\centerline{\sevenrm Faculty of Sciences, University of M\'alaga}
\centerline{\sevenrm Campus de Teatinos, S/N, 29080, M\'alaga, Spain}
\centerline{\sevenrm candido@apncs.cie.uma.es}
\bigskip
{\sevenrm In this note we study associative dialgebras proving that the most interesting such structures 
arise precisely when the algebra is not semiprime. In fact the presence of some \lq\lq perfection\rq\rq property
(simpleness, primitiveness, primeness or semiprimeness) imply that the dialgebra comes from an associative algebra
with both products $\iz$ and $\de$ identified. We also describe the class of zero-cubed algebras and apply its study
to that of dialgebras. Finally we describe two-dimensional associative dialgebras.}
\bigskip
\advance\sectno by 1
\centerline{\sl {\number\sectno}. Introduction}
\bigskip
The notion of Dialgebra is introduced and motivated by J. Loday in relation with Leibniz algebras.
A Leibniz algebra is a kind of  
 \lq\lq non-commutative Lie algebra\rq\rq. To be more precise,  
a Leibniz algebra is an algebra $A$ with product $[\ ,\ ]$ characterized by the so called Leibniz identity:
  $$[[x, y], z] = [[x, z], y] + [x, [y, z]].$$
 When the bracket $[\ ,\ ]$ is skew-symmetric, we get the definition of a Lie
algebra.  In the same way as
any associative algebra gives rise to a Lie algebra by antisymmetrization $[x, y]:= xy - yx$.
we can associate to any associative dialgebra a Leibniz algebra structure by a new kind of antisymmetrization.
The key point is to start  with two distinct operations for the product $xy$ and the product $yx$, so that
 the bracket is not necessarily skew-symmetric. So, we define
 an associative dialgebra  as a vector space $A$ provided with two associative operations $\iz$ and $\de$, 
called respectively left and right product, satisfying the identities:
$$x\iz (y \iz z) = x \iz (y \de z),$$
$$(x \de y) \iz z = x \de (y \iz z),$$
$$ (x \iz y) \de z = (x \de y) \de z.$$
Then one can easily check that the bracket $[x, y] := x \iz y - y \de x$ defines a Leibniz algebra. Therefore any associative
dialgebra gives rise to a Leibniz algebra. A classical example of dialgebra is constructed from a differential associative
algebra $(A, d)$ by defining $x\iz y:= x d(y)$ and $x\de y := d(x) y$. Then it is routinary to prove that
$(A, \iz, \de)$ is an associative dialgebra. There is a natural
dialgebra structure on the de Rham complex of a manifold.

Dialgebras appear in the literature in quite different contexts. So for instance in [\number\juana], dialgebras
are considered from the viewpoint of identities. However in [\number\baxter] they are studied from the Yang-Baxter
equation point of view. They can be related also to triple products as in [\number\sib] or to Leibniz algebras as
previously mentioned (see [\number\lodaydos]). 

\bigskip
\advance\sectno by 1
\centerline{\sl {\number\sectno}. Preliminaries on zero-cubed associative algebras}
\bigskip
One easy example of dialgebra arises when one considers an associative algebra $A$ with product denoted by
juxtaposition. Then we can define a dialgebra by writing $\iz:=0$ and defining $\de$ to be the product in $A$.
However the identities for an associative dialgebra imply that $A$ must be a nilpotent associative algebra of
index tree, that is, one must have $A^3=0$ if the identities for an associative dialgebra are to be satisfied.

So we start our study describing those associative algebras $A$ such that $A^3=0$. We call them {\sl zero-cubed} algebras.
Before this, we give a rather general example of associative algebra $A$ satisfying $A^3=0$. Consider two vector spaces
$Z$ and $X$ over the same field $F$ and suppose that there is a bilinear map $f\colon X\times X\to Z$. Define then in
the vector space $Z\oplus X$ the multiplication 
\advance\cf by 1\reone=\cf
$$
(z+x)(z'+x')=f(x,x'),\qquad z,z'\in Z, x,x'\in X.\eqno{(\number\cf)}$$ 
One can check immediately that $A:=Z\oplus X$ is an associative algebra such that $A^3=0$.
For any such algebra consider the annihilator $Z$
defined as the set of all elements $x\in A$ such that $xA=Ax=0$. In general, this is an ideal of the algebra. In our case
we have also $A^2\subset Z$. If we choose a complement $X$ for $Z$ in $A$ we have a direct sum of vector subspaces
$A=Z\oplus X$ and of course $X^2\subset A^2\subset Z$. Thus, there is a bilinear map $f\colon X\times X\to Z$ given by
$f(a,b):=ab$.
\medskip
\th{For any associative zero-cubed algebra $A$ (that is $A^3=0$) there are vector spaces $Z$ and $X$ over $F$ and a bilinear map
$f\colon X\times X\to Z$ such that $A=Z\oplus X$ with multiplication as in (\number\reone).
Furthermore take any two $F$-algebras $A=Z_A\oplus X_A$ and $B=Z_B\oplus X_B$ such that $A^3=0$, $B^3=0$. Denote
by $f_A\colon X_A\times X_A\to Z_A$ and $f_B\colon X_B\times X_B\to Z_B$ the bilinear maps involved in the products
of $A$ and $B$ as in (\number\reone), 
then $A\cong B$ if and only if there are isomorphisms $\alpha\colon Z_A\to Z_B$ and $\beta\colon X_A\to X_B$ such
that  $f_B(\beta(x),\beta(x'))=\alpha(f_A(x,x'))$.
}
\smallskip
Proof. It remains to prove the isomorphism condition. Thus suppose that $\omega\colon A\to B$ is an isomorphism.
Then since $Z_A$ is the annihilator of $A$ and similarly for $Z_B$, we have $\omega(Z_A)=Z_B$. Thus define
$\alpha=\omega\vert_{Z_A}$. Now consider the composition
$X_A\buildrel{i}\over\to A\buildrel{\omega}\over\to B\buildrel{\pi_{X_B}}\over\to X_B$ where
$i$ is the inclusion map and $\pi_{X_B}$ the canonical projection. Denote this composition by $\beta$.
Let us check that $\beta$ is an isomorphism. So take $x\in X_A$ and suppose that $\beta(x)=0$. This means
that $\omega(x)\in Z_B$ hence $x\in Z_A$ (but $x\in X_A$!). Thus $x=0$ proving that $\beta$ is a monomorphism.
To see that it is also an epimorphism take $y\in X_B$. Then $\omega^{-1}(y)=x+z$ with $x\in X_A$ and $z\in Z_A$.
Equivalently $\omega(x)=y-\omega(z)$ so that $\beta(x)=y$ hence $\beta$ is an epimorphism.
Finally if we take $x,x'\in X_A$ we know that $\omega(xx')=\omega(x)\omega(x')$. But
$xx'=f_A(x,x')\in Z_A$ and $\omega(xx')=\alpha(f_A(x,x'))$. On the other hand $\omega(x)=y+z_1$ and
$\omega(x')=y'+z_2$ where $y,y'\in X_B$, $z_1,z_2\in Z_B$. This implies that $\omega(x)\omega(x')=f_B(y,y')$
and of course $y=\beta(x)$, $y'=\beta(x')$ so that we have proved $f_B(\beta(x),\beta(x'))=\alpha(f_A(x,x'))$
as required. The theorem is proved.
\medskip
Consider the class $\cal C$ of all triples $(Z,X,f)$ where $Z$ and $X$ are vector spaces over $F$ and $f\colon X\times X\to Z$
a bilinear map. If $(Z_A,X_A,f_A)$ and $(Z_B,X_B,f_B)$ are two such triples we define the following relation:
$(Z_A,X_A,f_A)\equiv (Z_B,X_B,f_B)$ if there are isomorphisms $\alpha\colon Z_A\to Z_B$ and $\beta\colon X_A\to X_B$ such
that  $f_B(\beta(x),\beta(x'))=\alpha(f_A(x,x'))$. Then the isomorphy classes of zero-cubed associative $F$-algebras 
are in one-to-one correspondence with the equivalence classes of $\cal C$ modulo the relation $\equiv$.
\medskip
Recall that a semiprime (associative) algebra $A$ is one in which for any ideal $I$ of $A$
if $I^2=0$ then $I=0$. An easy observation is the following one:
\medskip
\co{Semiprime zero-cubed algebras do not exist.}
\newcount\etiq
\etiq=\coco
\smallskip
Proof. We are proving that if $A$ is a semiprime zero-cubed algebra then $A=0$.
If $A\ne 0$ is a zero-cubed algebra $A=Z\oplus X$ where $Z$ is the annihilator of $A$ and $X^2\subset Z$.
Since $A$ is semiprime $A^2\ne 0$ and since $Z^2=0$ we have $Z=0$. 
But then $A=X$ and $A^2=X^2\subset Z=0$
a contradiction.  The corollary is proved.
\medskip
\co{Prime zero-cubed algebras do not exist, hence simple zero-cubed algebras do not exist and primitive zero-cubed
algebras do not exist.}
\medskip
\newcount\etiquno
\etiquno=\coco
\newcount\labco

\co{Any zero-cubed two-dimensional associative $F$-algebra $A$ is either a trivial algebra (that is $A^2=0$) or isomorphic
to $F\times F$ with multiplication $(\lambda,\mu)(\lambda',\mu'):=(\mu\mu',0)$ for any $\lambda,\lambda',\mu,\mu'\in F$.}
\medskip\labco=\coco
Proof. If $Z=A$ then $A$ is a trivial algebra of zero product. If $Z=0$ since $A^2\subset Z=0$ we have again that
$A$ is a trivial algebra. If $\dim(Z)=1$ then $\dim(X)=1$ and taking generators $z\in Z$ and $x\in X$ we get a basis
with $x^2\in Fz$. So scaling $z$ if necessary we may suppose $x^2=z$ and the isomorphism follows directly.
\bigskip
\advance\sectno by 1
\centerline{\sl {\number\sectno}. Associative dialgebras}
\bigskip
For a fixed field $F$ an associative $F$-dialgebra is a vector space $A$ provided with two binary operations
$\iz,\de\colon A\times A\to A$ such that $(A,\iz)$ and $(A,\de)$ are associative algebras and 
 $$(x\dashv y)\dashv z\buildrel{(1)}\over{=}x\dashv (y\vdash z),\quad
(x\de y)\iz z\buildrel{(2)}\over{=}x\de(y\iz z),\quad
(x\iz y)\de z\buildrel{(3)}\over{=}x\de (y\de z).
$$
for any $x,y,z\in A$.
\medskip
There are several trivial ways to construct associative dialgebras from associative algebras. For instance if $A$
is an associative algebra with product denoted by juxtaposition, then we can define the dialgebra whose underlying
vector space agree with that of $A$ and whose products are $x\iz y:= x y=: x\de y$ for any $x,y\in A$. We shall call
this the {\sl associative dialgebra coming from the associative algebra} $A$. Often we shall rule out this algebra
since it reduces to an associative algebra. Another trivial way to construct associative dialgebras is by
considering one of the products to be zero. For instance if $\iz=0$ then $(A,\de)$ is an associative algebra and
furthermore, it is a zero-cubed algebra which has been studied in previous sections of this work. The same applies if $\de =0$.
\medskip
Some standard definitions of annihilators in an associative dialgebra are the following
$$\matrix{
\raniz(A):=\{x\in A\colon A\iz x=0\},& \laniz(A):=\{x\in A\colon x\iz A=0\}\cr
\rande(A):=\{x\in A\colon A\de x=0\},& \lande(A):=\{x\in A\colon x\de A=0\}}
$$
As usual  $\raniz(A),\laniz(A)$ are ideals of $(A,\iz)$ while 
$\rande(A)$ and $\lande(A)$ are ideals of $(A,\de)$.

An interesting definition on associative dialgebras is that of {\sl dual dialgebra}, $\op{A}$, of a given
associative dialgebra $A$.
\medskip
\defi{For a dialgebra $(A,\iz,\de)$ we denote by $\op{A}$ the dialgebra $(A,\iz',\de')$
whose underlying vector space agree with that of $A$ and whose new products are
$x\iz'y:=y\de x$ and $x\de' y:= y\iz x$ for any two elements $x,y\in A$. We will call $\op{A}$
the {\rm opposite algebra} of $A$.} 
\medskip
It is routinary to check that in case $A$ is an associative dialgebra then $\op{A}$ is 
also an associative dialgebra. Furthermore we have
$\raniz(A)=\lande(\op{A})$, $\rande(A)=\laniz(\op{A})$, $\laniz(A)=\rande(\op{A})$ and $\lande(A)=\raniz(\op{A})$.
In general we shall call {\sl duality} to the fact of interpreting a result (proved for an associative dialgebra $A$)
in the opposite algebra $\op{A}$.
\medskip
\pr{
In any associative dialgebra $A$ the annihilators $\raniz(A)$ and
$\lande(A)$ are ideals of $A$.}
\smallskip
Proof. 
By duality we only need to prove that $\raniz(A)$ is an ideal of $A$.
Let us denote $R=\raniz(A)$. We know that $R$ is an ideal of $(A,\iz)$ so we only
need to check that  $R\de A+A\de R\subset R$.
Let us  prove  $A\de R\subset R$ and $R\de A\subset R$. For the first we must see that 
$A\iz(A\de R)=0$. But applying (1) one has $A\iz(A\de R)=(A\iz A)\iz R\subset A\iz R=0$.
To see that $R\de A\subset R$ we must prove $A\iz (R\de A)=0$. But applying again (1) we have
$A\iz (R\de A)=(A\iz R)\iz A=0$ as we wanted to see. 
Thus we have proved $$\raniz(A)\triangleleft A.$$
Now that we know that $\raniz(A)$ and $\lande(A)$ are ideals of $A$ we can define the {\sl annihilator} of $A$
(denoted $\ann(A)$) as the ideal $$\ann(A):=\raniz(A)\cap\lande(A).$$ 
\rem{
In a dialgebra $A$ for any $y,z\in A$ we have $y\iz z-y\de z\in\raniz(A)$. Indeed, the identity (1) can
be written as $x\iz(y\iz z-y\de z)=0$ for any $x\in A$. Hence $y\iz z-y\de z\in\raniz(A)$.
Of course by duality we also have $y\iz z-y\de z\in\lande(A)$ hence $y\iz z-y\de z\in\ann(A)$ for
any $y,z\in A$.
}\kaka=1
\medskip
On the other hand we recall that a {\sl bar-unit} of a dialgebra $A$ is an element $e\in A$ such that
$x\iz e=x=e\de x$ for all $x\in A$.
\medskip
\th{If $A$ is a dialgebra such that $A\iz A=A$ (for instance if $A$ has a bar-unit)
 then the quotient dialgebra $A/\raniz(A)$ has zero right annihilator: 
$$\raniz(A/\raniz(A))=0.$$
Similar statement follow by duality for $\lande(A)$ and also $\ann(A/\ann(A))=0$.}\smallskip
Proof. Denote  by $\bar x=x+\raniz(A)$ the equivalence class of $x$, by $\bar A$ the quotient algebra
$\bar A:=A/\raniz(A)$ and suppose $\bar A\iz \bar x=\bar 0$. Then $A\iz x\in\raniz(A)$ so that $A\iz(A\iz x)=0$.
Then  $(A\iz A)\iz x=0$ hence $A\iz x=0$ implying $x\in\raniz(A)$ and $\bar x=\bar 0$. The remaining assertions
are proved in a similar way.
\medskip
\th{Let $A$ be an associative dialgebra, then:
\smallskip
\item{a)} If $\raniz(A)=A$ then $A$ with the operation $\de$ is zero-cubed algebra (that is $A$ can be described
saying that $\iz$ is the zero product and $(A,\de)$ is an in the previous section. 
\smallskip
\item{b)} If $\raniz(A)=0$ then $\iz$ and $\de$ agree and $A$ is the dialgebra coming from an associative algebra.
If $A$ has a bar unit $e$, then it is the dialgebra coming from an associative algebra with unit $e$.
}
\newcount\etiqth
\etiqth=\thco
\smallskip
Proof. If $\raniz(A)=A$ then $A\iz A=0$ so the $\iz$ product is null. Using (3) we get $A\de A\de A=0$ so that $(A,\de)$
is a zero-cubed algebra.
On the other hand if $\raniz(A)=0$ then by Remark {\number\kaka} we have $y\iz z-y\de z\in\raniz(A)=0$ hence both products coincide. 
So the dialgebra comes from an associative algebra $A$ with product
$\cdot$ defining both products $\iz$ and $\de$ to be $\cdot$. If the dialgebra $A$ has a bar unit then $A$ is a unital associative
algebra. The theorem is proved.
\medskip
As we have seen the maximal and minimal cases $\raniz(A)=A$ and $\raniz(A)=0$ are not interesting. The general situation will be 
the existence of a short exact sequence $$0\to\raniz(A)\to A\to A/\raniz(A)\to 0.$$
Thus, under the condition $A=A\iz A$ we can say that
$A$ is a extension of a associative algebra ($A/\raniz(A)$) by a zero-cubed algebra ($\raniz(A))$.
\medskip
\defi{A dialgebra $A$ with operations $\de$ and $\iz$ as usual is said to be $\de$-simple if $(A,\de)$ is simple and
is said to be $\iz$-simple if $(A,\iz)$ is it. The dialgebra $A$ will be called simple iff it is both
$\de$-simple and $\iz$-simple.}
\medskip
\th{For an associative dialgebra $A$ the following assertions are equivalent:
\item{i)} $A$ is $\de$-simple.
\item{ii)} $A$ is $\iz$-simple.
\item{iii)} $A$ is simple.
\item{iv)} $A$ is the dialgebra associated to a simple associative algebra.}
\smallskip
Proof.
First we observe that iv) implies i), ii) and iii). Next we prove that any of them implies iv). So suppose i),
then $\raniz(A)$ is $0$ or $A$. But if $\raniz(A)=A$ then $A\iz A=0$ hence $(A,\de)$ is a zero-cubed algebra by (3) which 
is impossible by Corollary \number\etiquno. Thus necessarily $\raniz(A)=0$ and iv) follows by Theorem \number\etiqth.
Suppose now ii), then again $\raniz(A)=0$ or $\raniz(A)=A$. In this second case $A\iz A=0$ contradicting the fact that
$A$ is $\iz$-simple. So necessarily $\raniz(A)=0$ and we conclude as before. Finally iii) implies i) which implies iv).
The theorem is proved.
\medskip
\defi{Given a property $P$ of associative algebras and an associative dialgebra $A$, we will say that 
$A$ has the $\de$-property $P$ if $(A,\de)$ has the property $P$. Similarly we define dialgebras with the $\iz$-property.
We will say that $A$ has the property $P$ iff $(A,\de)$ and $(A,\iz)$ have it.
Thus we can consider $\de$-prime, $\de$-semiprime or $\de$-primitive dialgebras (and dually $\iz$-prime, 
$\iz$-semiprime or $\iz$-primitive dialgebras). Also, with this meaning we consider semiprime, prime or primitive dialgebras.}
\medskip
\th{For an associative dialgebra $A$ and being $P$ the property of being semiprime, prime or primitive,
the following assertions are equivalent:
\item{i)} $A$ has the $\de$-property $P$.
\item{ii)} $A$ has the $\iz$-property $P$.
\item{iii)} $A$ has the property $P$.
\item{iv)} $A$ is the dialgebra associated to an associative algebra satisfying $P$.}
\smallskip
Proof.  We can argue as before that 
iv) implies i), ii) and iii). Next, we must take into account that primitive $\Rightarrow$ prime $\Rightarrow$ semiprime. 
To see that i) or ii) implies iv), if $(A,\de)$ has the property $P$,
since $\lande(A)^2=0$ we have $\lande(A)=0$ hence $A$ is the dialgebra coming from an associative algebra with
both products $\de$ and $\iz$ agreeing (this is the dual assertion of Theorem \number\etiqth.b). Of course iii) implies i) hence iii) implies iv).
Thus the theorem is proved.
\bigskip
\advance\sectno by 1
\centerline{\sl {\number\sectno}. Low dimensional associative dialgebras}
\bigskip

In this section we classify 2-dimensional associative dialgebras which do not come from an associative algebra (that
is both products $\de$ and $\iz$ are different). Thus $\ann(A)\ne 0$ implying $\raniz(A)\ne 0$ (see Theorem \number\etiqth).
Of course we exclude also the trivial case $A=\ann(A)$. Thus $\dim(\ann(A))=1$ so that the annihilator is generated by
an element (say $r$) and we can write $\ann(A)=F r$. Choosing any complement $s$ such that $A=Fr\oplus Fs$,
the multiplication tables of $A$ are
\bigskip

\hbox{
\vbox{\offinterlineskip
\tabskip=0pt
\halign{ 
\vrule height2.75ex depth1.25ex width 0.6pt #\tabskip=1em &
\hfil #\hfil &\vrule # &  #\hfil &\vrule # &
\hfil #\hfil &#\vrule width 0.6pt \tabskip=0pt\cr
\noalign{\hrule height 0.6pt}
& $\iz$ && $r$ &&  $s$ & \cr
\noalign{\hrule}
& $r$ && $0$ && $x_1 r$ &\cr
& $s$ && $0$ && $x_2 r+x_3 s$ &\cr
\noalign{\hrule height 0.6pt}
}} 
\hskip 3cm
\vbox{\offinterlineskip
\tabskip=0pt
\halign{ 
\vrule height2.75ex depth1.25ex width 0.6pt #\tabskip=1em &
\hfil #\hfil &\vrule # &  #\hfil &\vrule # &
\hfil #\hfil &#\vrule width 0.6pt \tabskip=0pt\cr
\noalign{\hrule height 0.6pt}
& $\de$ && $r$ &&  $s$ & \cr
\noalign{\hrule}
& $r$ && $0$ && $0$ &\cr
& $s$ && $x_4 r$ && $x_5 r+x_6 s$ &\cr
\noalign{\hrule height 0.6pt}
}}
}
\medskip
Imposing the associativity conditions to the products $\iz$, $\de$ as well as the conditions (1), (2) and (3),
we get the following equations:
$$0=x_1x_2=x_4x_5=x_1(x_1-x_3)=x_4(x_3-x_4)=x_1(x_1-x_6)=x_2(x_1+x_3-x_6)=$$
$$x_1x_5+x_2x_6-x_2x_4-x_3x_5=x_3(x_3-x_6)=x_4(x_4-x_6)=x_6(x_3-x_6)=x_5(x_3-x_4-x_6)$$ 
and we analyze the different possibilities:
\medskip
\item{1)} $x_1=x_2=0$. Then if $x_3=0$ the product $\iz$ is null and $(A,\de)$ is a zero-cubed
algebra described in Corollary \number\labco. So we suppose $x_3\ne 0$ which implies $x_3=x_6$.
The identity $x_1x_5+x_2x_6-x_2x_4-x_3x_5=0$ implies  $x_5=0$
The previous list of identities reduces then to $0=x_4(x_3-x_4)$, $x_3=x_6$, $x_5=0$.
We have two subcases:
\item{1.i)} $x_4=0$, $x_3=x_6\ne 0$, $x_5=0$. Then $(A,\iz)$ is given by $0=r\iz r=r\iz s=s\iz r$ and $s\iz s=x_3s$.
Defining $s':=x_3^{-1}s$ the multiplication table relative to the basis $\{r,s'\}$ is $r\iz r=r\iz s'=s'\iz r=0$ and
$s'\iz s'=x_3^{-2}x_3s=x_3^{-1}s=s'$. Moreover $r\de r=r\de s'=s'\de r=0$ and $s'\de s'=s'$ as for the product $\iz$.
So in this case both operations $\iz$ and $\de$ coincide and $A$ is the associative dialgebra coming from the bidimensional
associative algebra with multiplication table: \medskip
\centerline{\vbox{\offinterlineskip
\tabskip=0pt
\halign{ 
\vrule height2.75ex depth1.25ex width 0.6pt #\tabskip=1em &
\hfil #\hfil &\vrule # &  #\hfil &\vrule # &
\hfil #\hfil &#\vrule width 0.6pt \tabskip=0pt\cr
\noalign{\hrule height 0.6pt}
& $\cdot$ && $r$ &&  $s$ & \cr
\noalign{\hrule}
& $r$ && $0$ && $0$ &\cr
& $s$ && $ 0$ && $s$ &\cr
\noalign{\hrule height 0.6pt}
}}}
\medskip
\item{1.ii)} $x_3=x_4=x_6\ne 0$, $x_5=0$. Defining again $s'=x_3^{-1}s$ we get as before 
$r\iz r=r\iz s'=s'\iz r=0$ and $s'\iz s'=s'$. Moreover the product $\de$ satisfies $r\de s'=0$, $s'\de r=r$ and
$s'\de s'=s'$. Hence after scaling $s$ if necessary we get the dialgebra with operations:
\medskip
\hbox{\hskip 1cm
\vbox{\offinterlineskip
\tabskip=0pt
\halign{ 
\vrule height2.75ex depth1.25ex width 0.6pt #\tabskip=1em &
\hfil #\hfil &\vrule # &  #\hfil &\vrule # &
\hfil #\hfil &#\vrule width 0.6pt \tabskip=0pt\cr
\noalign{\hrule height 0.6pt}
& $\iz$ && $r$ &&  $s$ & \cr
\noalign{\hrule}
& $r$ && $0$ && $0$ &\cr
& $s$ && $0$ && $s$ &\cr
\noalign{\hrule height 0.6pt}
}}
 \hskip 5cm
\vbox{\offinterlineskip
\tabskip=0pt
\halign{ 
\vrule height2.75ex depth1.25ex width 0.6pt #\tabskip=1em &
\hfil #\hfil &\vrule # &  #\hfil &\vrule # &
\hfil #\hfil &#\vrule width 0.6pt \tabskip=0pt\cr
\noalign{\hrule height 0.6pt}
& $\de$ && $r$ &&  $s$ & \cr
\noalign{\hrule}
& $r$ && $0$ && $0$ &\cr
& $s$ && $r$ && $s$ &\cr
\noalign{\hrule height 0.6pt}
}}
}
This algebra will be denoted by $\hbox{\bf I}$.
\medskip
\item{2)} $x_1=0$, $x_2\ne 0$. The list of identities reduces to
$0=x_4x_5=x_4(x_3-x_4)$, $x_3=x_6$,
$x_2x_3-x_2x_4-x_3x_5=0$. We consider the following subcases:
\item{2.i)} $x_1=0$, $x_2\ne 0$, $x_4=0$. Then we have $x_3=x_6$ and $x_3(x_2-x_5)=0$ and there
are again two possibilities:
\parindent=1cm
\item{2.i.a)} $x_1=0$, $x_2\ne 0$, $x_4=0$, $x_3=x_6=0$. After scaling $r$ if necessary we get the family
of dialgebras with multiplications given by
\medskip
\hbox{\hskip 1cm
\vbox{\offinterlineskip
\tabskip=0pt
\halign{ 
\vrule height2.75ex depth1.25ex width 0.6pt #\tabskip=1em &
\hfil #\hfil &\vrule # &  #\hfil &\vrule # &
\hfil #\hfil &#\vrule width 0.6pt \tabskip=0pt\cr
\noalign{\hrule height 0.6pt}
& $\iz$ && $r$ &&  $s$ & \cr
\noalign{\hrule}
& $r$ && $0$ && $0$ &\cr
& $s$ && $0$ && $r$ &\cr
\noalign{\hrule height 0.6pt}
}}
 \hskip 5cm
\vbox{\offinterlineskip
\tabskip=0pt
\halign{ 
\vrule height2.75ex depth1.25ex width 0.6pt #\tabskip=1em &
\hfil #\hfil &\vrule # &  #\hfil &\vrule # &
\hfil #\hfil &#\vrule width 0.6pt \tabskip=0pt\cr
\noalign{\hrule height 0.6pt}
& $\de$ && $r$ &&  $s$ & \cr
\noalign{\hrule}
& $r$ && $0$ && $0$ &\cr
& $s$ && $0$ && $ k r$ &\cr
\noalign{\hrule height 0.6pt}
}}
}
We shall denote this dialgebras by $\hbox{\bf II}_k$.
\medskip
\item{2.i.b)} $x_1=0$, $x_2\ne 0$, $x_4=0$, $x_3=x_6\ne 0$, $x_2=x_5$. In this case again the product coincide
and the algebra comes from an associative algebra.
\medskip
\item{2.ii)} $x_1=0$, $x_2\ne 0$, $x_4\ne 0$. This implies $x_5=0$, $x_3=x_4=x_6\ne 0$. Then replacing $s$ with
$s':=(x_2/x_3^2)r+x_3^{-1}s$ we get $r\iz r=r\iz s'=s'\iz r=0$ and $s'\iz s'=s'$. On the other hand
$r\de r=r\de s'=0$ and $s'\de r=r$, $s'\de s'=s'$. But this algebra is isomorphic to $\hbox{\bf I}$.
\medskip
\item{3)} $x_1\ne 0$, $x_2=x_4=0$. The list of identities is then
$x_1=x_3=x_6$. Thus necessarily $x_6\ne 0$ and we can define $s':=x_6^{-2}x_5r+x_6^{-1}s$. Then replacing $s$ with $s'$ if
necessary we get the dialgebra with multiplications:
\medskip
\hbox{\hskip 1cm
\vbox{\offinterlineskip
\tabskip=0pt
\halign{ 
\vrule height2.75ex depth1.25ex width 0.6pt #\tabskip=1em &
\hfil #\hfil &\vrule # &  #\hfil &\vrule # &
\hfil #\hfil &#\vrule width 0.6pt \tabskip=0pt\cr
\noalign{\hrule height 0.6pt}
& $\iz$ && $r$ &&  $s$ & \cr
\noalign{\hrule}
& $r$ && $0$ && $r$ &\cr
& $s$ && $0$ && $s$ &\cr
\noalign{\hrule height 0.6pt}
}}
 \hskip 5cm
\vbox{\offinterlineskip
\tabskip=0pt
\halign{ 
\vrule height2.75ex depth1.25ex width 0.6pt #\tabskip=1em &
\hfil #\hfil &\vrule # &  #\hfil &\vrule # &
\hfil #\hfil &#\vrule width 0.6pt \tabskip=0pt\cr
\noalign{\hrule height 0.6pt}
& $\de$ && $r$ &&  $s$ & \cr
\noalign{\hrule}
& $r$ && $0$ && $0$ &\cr
& $s$ && $0$ && $s$ &\cr
\noalign{\hrule height 0.6pt}
}}
} 
\medskip
We denote this algebra by $\hbox{\bf III}$.
\medskip
\item{4)} $x_1\ne 0$, $x_2=0$, $x_4\ne 0$, $x_5=0$. The list of identities is then
$x_1=x_3=x_4=x_6$. Replacing $s$ with $s':=x_1^{-1}s$ if necessary we find the dialgebra with multiplications
\medskip
\hbox{\hskip 1cm
\vbox{\offinterlineskip
\tabskip=0pt
\halign{ 
\vrule height2.75ex depth1.25ex width 0.6pt #\tabskip=1em &
\hfil #\hfil &\vrule # &  #\hfil &\vrule # &
\hfil #\hfil &#\vrule width 0.6pt \tabskip=0pt\cr
\noalign{\hrule height 0.6pt}
& $\iz$ && $r$ &&  $s$ & \cr
\noalign{\hrule}
& $r$ && $0$ && $r$ &\cr
& $s$ && $0$ && $s$ &\cr
\noalign{\hrule height 0.6pt}
}}
 \hskip 5cm
\vbox{\offinterlineskip
\tabskip=0pt
\halign{ 
\vrule height2.75ex depth1.25ex width 0.6pt #\tabskip=1em &
\hfil #\hfil &\vrule # &  #\hfil &\vrule # &
\hfil #\hfil &#\vrule width 0.6pt \tabskip=0pt\cr
\noalign{\hrule height 0.6pt}
& $\de$ && $r$ &&  $s$ & \cr
\noalign{\hrule}
& $r$ && $0$ && $0$ &\cr
& $s$ && $r$ && $s$ &\cr
\noalign{\hrule height 0.6pt}
}}
} 
\medskip
We denote this algebra by $\hbox{\bf IV}$.
\medskip

Now we check that the algebras $\hbox{\bf I}$, $\hbox{\bf II}_k$, $\hbox{\bf III}$ and $\hbox{\bf IV}$
are all non-isomorphic. Thus denote by $A^{2\ \iz}:=A\iz A$ and similarly $A^{2\ \de}:=A\de A$. Then we compute
the dimension of these spaces in the following table: 
\medskip
\centerline{
\vbox{\offinterlineskip
\tabskip=0pt
\halign{ 
\vrule height2.75ex depth1.25ex width 0.6pt #\tabskip=1em &
\hfil #\hfil &\vrule # &  #\hfil &\vrule # &
\hfil #\hfil &#\vrule width 0.6pt \tabskip=0pt\cr
\noalign{\hrule height 0.6pt}
& $A$ && $\dim(A^{2\ \iz})$ && $\dim(A^{2\ \de})$  & \cr
\noalign{\hrule}
& $\hbox{\bf I}$ && $1$ && $2$ &\cr
& $\hbox{\bf II}_k,\ (k\ne 0)$ && $1$ && $1$ &\cr
& $\hbox{\bf III}$ && $2$ && $1$ &\cr
& $\hbox{\bf IV}$ && $2$ && $2$ &\cr
\noalign{\hrule height 0.6pt}
}}}
\medskip
\noindent which proves that the algebras are not isomorphic except possibly for the family of algebras $\hbox{\bf II}_k$
($k\ne 0$).
Let us now check that for these algebras, $\hbox{\bf II}_k\cong \hbox{\bf II}_{k'}$ if and only if $k=k'$.
For this we shall need to compute the group of automorphisms of the two-dimensional algebra $A$ with basis
$\{r,s\}$ such that $r^2=rs=sr=0$ and $s^2=r$. If $f\in\hbox{aut}(A)$ then necessarily $f(r)=k_1 r$ for a nonzero 
scalar $k_1\in F$. Then, since $s^2=r$ we may write $f(s)=\alpha r+\beta s$ for some scalars $\alpha,\beta\in F$ and
we have $(\alpha r+\beta s)^2=k_1 r$ which implies $\beta^2 r=k_1 r$ so that $\beta^2=k_1$. Thus
$$\Bigg\{
\matrix{ 
f(r)=\beta^2 r \cr
f(s)=\alpha r+\beta s.}$$
As a corollary $\hbox{aut}(A)$ is isomorphic to the subgroup of $\hbox{GL}_2(F)$ of all matrices of the form
$$\pmatrix{\beta^2 & 0\cr \alpha & \beta}$$
where $\alpha\in F$, $\beta\in F^\times$. Thus if we have an isomorphism $f\colon\hbox{\bf II}_k\cong\hbox{\bf II}_{k'}$
then we have an isomorphism $f\colon (\hbox{\bf II}_k,\iz)\cong (\hbox{\bf II}_{k'},\iz)$ which implies that $f$
must have the canonical form given above. But on the other hand $f$ is also an isomorphism
$(\hbox{\bf II}_k,\de)\cong (\hbox{\bf II}_{k'},\de)$ which implies $f(s\de s)=f(s)\de f(s)$ that is
$kf(r)=(\alpha r+\beta s)\de (\alpha r+\beta s)$ or equivalently
$k\beta^2 r=\beta^2 k' r$ which implies $k=k'$. Summarizing all the information in this section we claim:
\medskip
\th{Let $A$ be a two-dimensional associative dialgebra over a field $F$. Then we have only one of the following possibilities
for $A$:\smallskip
\item{i)} $A$ is the dialgebra coming from an associative algebra.\smallskip
\item{ii)} One (and only one) of the products $\iz$ or $\de$ is null. Then $A$ is a zero-cubed algebra with the nonzero product 
and has been described in Corollary \number\labco.\smallskip
\item{iii)} $A$ is isomorphic to one of the algebras $\hbox{\bf I}$, $\hbox{\bf II}_k$ ($k\ne 0$), $\hbox{\bf III}$ or $\hbox{\bf IV}$
previously described. Any two of these algebras are not isomorphic.
}
\medskip
\centerline{\bf Acknowledgements}
\medskip
The author would like to thank J. S\'anchez Ortega for lectures given at the {\sl M\'alaga seminar
of the research project FQM 336}, where the notion of dialgebra was introduced 
and motivated to this researcher. Thanks are also given to D. Mart\'{\i}n Barquero for correcting
an earlier version of the work, which improved remarkably the paper. 
\medskip
\centerline{\bf References}
\medskip
\parindent=0cm

[\number\juana] 
M. R. Bremner and J. S\'anchez Ortega. {\it The partially alternating ternary sum in an
associative dialgebra}. J. Phys. A: Math. Theor. 43 (2010) 455215 
doi:10.1088/1751-8113/43/45/455215.
\smallskip

[\number\baxter] 
Ra\'ul Felipe, Nancy L\'opez-Reyes, Fausto Ongaya and  Ra\'ul Vel\'asquez. {\it Yang-Baxter equations on matrix dialgebras with 
a bar unit}. Linear Algebra and its Applications
Volume 403, 1, 2005, p. 31-44.\smallskip

[\number\koles]
P. Kolesnikov. {\it Varieties of dialgebras and conformal algebras}. Siberian Math. J. 49, 257-272, 2008.
\smallskip

[\number\lodayuno]
J. Loday. {\it Alg\`ebres ayant deux op\`erations associatives (dig\`ebres)}. C. R. Acad. Sci. Paris S\`er ́I. Math. 321, 141-146, 1995.
\smallskip

[\number\lodaydos]
J. Loday. {\it Une version non commutative des alg\`ebres de Lie: les alg\`ebres de Leibniz}.  Enseign. Math. 39, 269-293, 1993.
\smallskip

[\number\operads]
J.-L. Loday, A. Frabetti, F. Chapoton and F. Goichot. {\it Dialgebras and Related Operads}. Lecture Notes in Mathematics 1763.
Springer-Verlag. 2001.\smallskip

[\number\sib]
A. Pozhidaev.  {\it Dialgebras and related triple systems}. Siberian Math. J. 49, 696-708, 2008.

\end